\documentclass[12pt]{article}

\usepackage{amssymb}

\newtheorem{theorem}{Theorem}[section]

\newtheorem{assumption}[theorem]{Assumption}

\newenvironment{proof}{\noindent {\bf
Proof}.\ }{\proofbox\par\smallskip\par}

\def\numberlikeadb{\global\def\theequation{\thesection.\arabic{equation}}}
\numberlikeadb

\newcommand{\halmos}{\rule{1ex}{1.4ex}}
\newcommand{\proofbox}{\hspace*{\fill}\mbox{$\halmos$}}

\newcommand{\pr}{\mbox{\bf P}}

\newcommand{\reals}{{R}}
\newcommand{\ints}{\mathbb{Z}}

\def\ZZ{\ints}

\newcommand{\eqa}{\begin{eqnarray}}
\newcommand{\ena}{\end{eqnarray}}
\newcommand{\eq}{\begin{equation}}
\newcommand{\en}{\end{equation}}
\newcommand{\eqs}{\begin{eqnarray*}}
\newcommand{\ens}{\end{eqnarray*}}

\def\a{\alpha}
\def\b{\beta}
\def\g{\gamma}

\def\d{\delta}

\def\e{\varepsilon}
\def\h{\eta}
\def\z{\zeta}

\def\k{\kappa}
\def\m{\mu}
\def\n{\nu}

\def\r{\rho}
\def\s{\sigma}

\def\nin{\noindent}
\def\bsk{\bigskip}
\def\msk{\medskip}

\def\Blb{\left\{}
\def\Brb{\right\}}

\def\non{\nonumber}
\def\integ{{\mathbb Z}}

\def\Eq{\ =\ }
\def\DEq{\ :=\ }
\def\Le{\ \le\ }

\def\sjo{\sum_{j\ge0}}

\def\sji{\sum_{j\ge1}}
\def\stT{\sup_{0\le t\le T}}
\def\sio{\sum_{i\ge0}}

\def\slo{\sum_{l\ge0}}

\def\Ref#1{{\rm (\ref{#1})}}
\def\vol#1{{\bf #1},}
\def\bp{\begin{proof}}
\def\ep{\end{proof}}

\def\bone{{\bf 1}}

\def\Bl{\left(}
\def\Br{\right)}

\def\half{{\textstyle{\frac12}}}

\def\xx{{\cal X}}

\def\ignore#1{}

\def\sjni{\sum_{j\neq i}}

\def\tW{{\widetilde W}}

\def\Z{\ints}
\def\Zp{\Z_+}

\def\R{\reals}
\def\rr{{\mathcal R}}
\def\jj{{\mathcal J}}
\def\TT{{\mathcal T}}
\def\sjj{\sum_{J\in \jj}}
\def\x{\xi}
\def\smo{\sum_{m\ge0}}
\def\sko{\sum_{k\ge0}}

\def\nm#1{{\|#1\|_\m}}
\def\uii{^{(i)}}

\def\Def{\DEq}

\def\tm{{\widetilde m}}

\def\rmax{r_{\rm max}}

\def\ui{^{(1)}}
\def\uii{^{(i)}}
\def\ut{^{(2)}}
\def\uh{^{(3)}}

\def\halfh{{\textstyle{\frac14}}}

\def\Blm{\left|}
\def\Brm{\right|}

\def\tE{{\widetilde E}}
\def\bx{{\bar x}}

\begin{document}

\title{Central limit approximations for Markov population
processes with countably many types}
\author{
A. D. Barbour\footnote{Angewandte Mathematik, Universit\"at Z\"urich,
Winterthurertrasse 190, CH-8057 Z\"URICH;
ADB was supported in part by Schweizerischer Nationalfonds Projekt Nr.\
20--107935/1 and by Australian Research Council Grants Nos DP120102728 and DP120102398.\msk}
\ and
M. J. Luczak\footnote{University of Sheffield;
MJL was supported in part by an EPSRC Leadership Fellowship No EP/J004022/1 and by Australian Research Council 
Grant No DP120102398.
}\\
Universit\"at Z\"urich and University of Sheffield
}

\date{}
\maketitle

\begin{abstract}
When modelling metapopulation dynamics, the influence of a single patch
on the metapopulation depends on the number of individuals in the patch.
Since there is usually no obvious natural upper limit on the number
of individuals in a patch, this leads to
systems in which there are countably infinitely many possible types of
entity.  Analogous considerations apply in the transmission of
parasitic diseases.  In this paper, we prove central limit theorems for
quite general systems of this kind, together with bounds on the rate
of convergence in an appropriately chosen weighted $\ell_1$ norm.
\end{abstract}

{\it Running head:}  A central limit approximation

\section{Introduction}\label{introduction}
\setcounter{equation}{0}

Metapopulations, introduced by Levins~(1969), are used to describe the evolution 
of the population of a species in a fragmented habitat. The metapopulation consists 
of a number of distinct patches,
together with (a summary of) the population present in each patch, and
its development over time is governed by specified within and between patch dynamics. 
In the Markovian structured mean-field metapopulation model of Arrigoni~(2003), 
the state of the system consists simply of the numbers of individuals in each patch.
Individuals reproduce within patches and migrate between patches, and each
patch is subject to random catastrophes, which reduce its population to zero.
Letting~$N$ be the total number of patches, thought of as being large, and letting
$X_t^{N,i}$ denote the number of patches with~$i$ individuals at time~$t$,
the transitions out of state~$X \in \ZZ_+^{\ZZ_+}$ to states $X+J$ in her model are 
as follows:
$$
\begin{array}{rlll}
  J &=\ e^{(i-1)} - e\uii &\mbox{at rate}\ Ni x^i(d_i + \g(1-\r)),\quad &i\ge2;\\[1ex]
  J &=\ e^{(0)} - e\ui &\mbox{at rate}\ N x^1(d_1 + \g(1-\r) + \k);\\[1ex]
  J &=\ e^{(i+1)} - e\uii &\mbox{at rate}\ Ni  x^i b_i,
       \quad &i \ge 1;\\[1ex]
  J &=\ e^{(0)} - e\uii &\mbox{at rate}\ N x^i \k, \quad &i\ge2;\\[1ex]
  J &=\ e^{(k+1)} - e^{(k)} + e^{(i-1)} - e\uii &\mbox{at rate}\ N i x^i x^k \r\g,
     \qquad\qquad k\ge0,\ &i\ge1;
\end{array}
$$
here, $x := N^{-1}X$.
The total number~$N := \sjo X^{N,j}$ of patches remains constant throughout, and the number
of individuals in any one patch changes by at most one at each transition.
The {\it per capita\/} death and birth rates $(d_i), (b_i)$ within each patch are
allowed to vary with the current population size~$i$, but in the same way in all patches;
they would usually be chosen to correspond to one of the traditional single species
demographic models.
The {\it per capita\/} migration rate~$\g$ is also the same for all individuals, as
is the probability~$\r$ that a migration is successful, and a successful migrant chooses
its new patch uniformly at random.  Each patch is independently subject to catastrophes
at the same rate~$\k$.

If there were an absolute upper limit for the number of individuals in each patch, 
the model would be a finite dimensional Markov population process.  The behaviour
of these finite dimensional models
can be approximated using the methods pioneered by Kurtz (1970, 1971), who was able to
establish a law of large numbers approximation, in the form of a system of ordinary
differential equations, and a corresponding diffusion approximation.  However, there
are no upper limits on population number in the usual single population models, and 
it is the stochastic evolution according to the rules of the model that dictates 
the region in which population
numbers typically lie.  Thus it seems unnatural to introduce an {\it a
priori\/} upper limit in the system above, just because more than one population 
is being modelled.  The same 
considerations surface in a number of other population models,
including the epidemic models of Luchsinger~(2001a,b) and Kretzschmar~(1993), and the
model of cell behaviour as a function of the copy number of a particular gene in
Kimmel \& Axelrod~(2002, Chapter~7). Instead, it makes sense to consider Markov
population processes with a countably infinite number of dimensions as models in 
their own right.

A law of large numbers in a general setting of this kind was first established by
Eibeck \& Wagner~(2003).  Under appropriate conditions, Barbour \& Luczak~(2008, 2011)
strengthened the law of large numbers by providing an error bound, in a weighted
$\ell_1$ norm, that is close to optimal order in~$N$.
In this paper, these latter results are complemented by a central limit
approximation, together with a corresponding error estimate.

Our general setting, as in Barbour \& Luczak~(2011) [BL], is that of families of 
Markov population processes $X^N := (X_t^N,\,t\ge0)$, $N\ge1$, taking values in the
countable space $\xx_+ := \{X \in \Zp^{\Zp};\, \smo X^m < \infty\}$.
The component $X_t^{N,j}$ of~$X_t^N$ represents the number of individuals of type~$j$
that are present at time~$t$, and there are countably many types possible; however, 
at any given time, there are only finitely many individuals in the system.  
The process evolves as a Markov process with state--dependent transitions
\eq\label{1.0}
    X \ \to\ X + J \quad\mbox{at rate}\quad N\a_J(N^{-1}X),\qquad X \in \xx_+,\ J \in \jj,
\en
where each jump is of bounded influence, in the sense that
\eq\label{finite-jumps}
   \jj \subset \Bigl\{X \in \Z^{\Zp};\, \smo |X^m| \le J_* < \infty \Bigr\},
     \quad\mbox{for some fixed}\quad J_* < \infty, 
\en
so that the number of individuals affected at each transition is uniformly bounded.  Density
dependence is reflected in the fact that the arguments of the functions~$\a_J$
are counts normalised by the `typical size'~$N$.
Writing $\rr := \R_+^{\Zp}$, the
functions~$\a_J\colon\, \rr  \to \R_+$ are assumed to satisfy
\eq\label{finite-rates}
    \sjj \a_J(\x) \ < \ \infty,\qquad \x \in \rr_0,
\en
where $\rr_0 := \{\x \in \rr\colon\, \x_i = 0 {\mbox{ for all but finitely many }} i\}$;
this assumption implies that the processes~$X^N$ are indeed pure jump processes, at least
for some non-zero length of time.  To prevent
the paths leaving~$\xx_+$, we also assume that $J^l \ge -1$ for each~$l$, and
that $\a_J(\x) = 0$ if $\x^l = 0$ for any $J\in\jj$ such that $J^l = -1$.

In the finite dimensional case, the law of large numbers is expressed in terms of the 
system of {\it deterministic equations\/}
\eq\label{determ}
   \frac{d\x}{dt} \Eq \sjj J\a_J(\x).
\en
In~[BL], it is assumed that
\eq\label{F-assn}
   \sjj J\a_J(\x) \Eq A\x + F(\x),
\en
where~$A$ is a constant $\Zp\times\Zp$ matrix, and~\Ref{determ} is then treated 
as a perturbed linear system  (Pazy~1983, Chapter~6). 
Under suitable assumptions on $A$, there exists a measure~$\m$ on~$\Zp$, 
defining a weighted $\ell_1$ norm $\| \cdot \|_{\mu}$ on ${\mathcal R}$,  
and a strongly $\| \cdot \|_{\mu}$--continuous semigroup $\{R(t),\,t\ge0\}$ 
of transition matrices having pointwise derivative $R'(0) = A$.  
If $F$ is locally $\| \cdot \|_{\mu}$--Lipschitz, the solution~$x$ of the integral 
equation
\eq\label{x-eqn}
   x_t \Eq R(t)x_0 + \int_0^t R(t-s) F(x_s)\,ds\,,
\en
for $\nm{x_0} < \infty$, replaces that of~\Ref{determ}
as an approximation to~$x^N := N^{-1}X^N$.

Under suitable conditions, it is shown in~[BL, Theorem~4.7] that
\[
  \stT \nm{x_t^N - x_t} \Eq O(N^{-1/2}\sqrt{\log N}),
\]
except on an event of probability of order $O(N^{-1}\log N)$, provided that
$\nm{x_0^N - x_0} = O(N^{-1/2}\sqrt{\log N})$.   The conditions under which
this approximation holds can be divided into three categories: growth conditions
on the transition rates, so that the {\it a priori\/} bounds, which have
the character of moment bounds, can be established; conditions on the matrix~$A$,
sufficient to limit the growth of the semigroup~$R$, and (together with the 
properties of~$F$) to determine the weights
defining the metric in which the approximation is to be carried out; and conditions
on the initial state of the system.  The conditions are described in the next section.
They are all needed in the
current paper, too, in which we investigate the difference $x_t^N - x_t$ in greater
detail.  

Our main result, Theorem~\ref{diffusion-theorem}, shows that, under some extra conditions,
it is possible to construct a diffusion process~$Y$ on the same probability space
as~$X^N$ in such a way that
\eq\label{CLT-1}
      \stT \nm{N^{1/2}(x_t^N - x_t) - Y_t} \Eq O(N^{-b_1}), 
\en
except on an event of probability of order~$O(N^{-b_2})$, for specific values of $b_1$
and~$b_2$.  With the best possible control of moments, as for the model of Arrigoni~(2003)
mentioned above, one can take any $b_1 < 1/4$ and any $b_2 < 1$, provided that the
initial conditions are appropriately chosen.
The process~$Y$ can be interpreted as the infinite dimensional
analogue of the diffusion approximation in Kurtz~(1971), satisfying the formal stochastic
differential equation
\eq\label{diffusion-eq1}
   dY_t \Eq \{A + DF(x_t)\}Y_t\,dt + dW_t.
\en
Here, $dW$ is a time--inhomogeneous white noise process with infinitesimal covariance 
matrix $\s^2(t) := \sjj JJ^T \a_J(x_t)$, and~$Y$ has time--inhomogeneous linear drift
with coefficient matrix $A + DF(x_t)$.  In particular, if~$\bx$ is an equilibrium of the
deterministic equations, satisfying $A\bx + F(\bx) = 0$, then~$Y$ is an infinite dimensional
Ornstein--Uhlenbeck process, with constant drift coefficient matrix $A + DF(\bx)$ and
infinitesimal covariance matrix $\sjj JJ^T \a_J(\bx)$.

\subsection*{Basic approach}\label{basic}
The structure of the argument is as follows. 
It is shown in [BL, (4.8)] that, under suitable conditions, the process~$x^N$ satisfies an equation
very similar to~\Ref{x-eqn}:
\eq\label{xN-eqn}
   x_t^N \Eq R(t)x_0^N + \int_0^t R(t-s) F(x_s^N)\,ds + \tm_t^N,
\en
where
\eq\label{mN-tilde-def}
   \tm_t^N \Def \int_0^t R(t-s)\,dm_s^N \Eq m_t^N + \int_0^t R(t-s)Am_s^N\,ds,
\en
and
\eq\label{mN-def}
    m_t^N \Def x_t^N - x_0^N - \int_0^t \{Ax_s^N + F(x_s^N)\}\,ds
\en
is a local martingale.  
Taking the difference between \Ref{xN-eqn} and~\Ref{x-eqn}, and multiplying
by~$\sqrt N$, gives
\eqa
  U_t^N &=& R(t) U_0^N + \int_0^t R(t-s)\,DF(x_s) [U_s^N]\,ds
     + \h_t^N + N^{1/2}\tm_t^N, \label{main-eqn}
\ena
with $U_t^N := N^{1/2}\{x_t^N - x_t\}$ and
\eq\label{eta-def}
   \h_t^N \Def N^{1/2}\int_0^t R(t-s)\{F(x_s^N) - F(x_s) - DF(x_s)[N^{-1/2}U_s^N]\}\,ds.
\en

Starting with this representation of~$U^N$, the first step is to show that~$\h^N$ is
uniformly small with high probability, so that the randomness in~$U^N$ is driven
principally by the process $N^{1/2}\tm^N$. This quantity is in turn determined,
through~\Ref{mN-tilde-def}, by the local martingale $N^{1/2}m^N$.  The next step
is to show that $N^{1/2}m^N$ is close to a diffusion~$W$, formally expressible
as
\eq\label{W-def}
    W_t \Def \sjj JW_J(A_J(t)),
\en
where the $\{W_J,\,J\in\jj\}$ are independent standard Brownian motions, and
$A_J(t) := \int_0^t \a_J(x_s)\,ds$: this is the diffusion appearing in~\Ref{diffusion-eq1}. 
Analogously to~\Ref{mN-tilde-def}, we then show that we can define a process~$\tW$ such that
\eq\label{tW-def}
   \tW_t \Def W_t + \int_0^t R(t-s)A W_s\,ds,
\en
and that $\tW$ is close to $N^{1/2}\tm^N$.  Finally, returning to~\Ref{main-eqn},
we show that $U^N$ is close to the solution~$Y$ to the analogous equation
\eq\label{Y-eqn}
  Y_t \Eq R(t) Y_0 + \int_0^t R(t-s)\,DF(x_s) [Y_s]\,ds  + \tW_t,
\en
which in turn can be shown to exist and be unique.
The random process solving~\Ref{Y-eqn} at first sight seems rather mysterious.
However, partial integration represents~$\tW_t$ as $\int_0^t R(t-s)\,dW_s$, and so
the expression for~$Y$ can indeed be interpreted as the variation of constants representation
of the solution to the formal stochastic differential equation~\Ref{diffusion-eq1}.

In the remaining sections, this programme is carried out in detail.  Section~\ref{assumptions}
is concerned with specifying the conditions under which the main theorem is true, and with
recalling the results from~[BL] that are needed here.  In the subsequent sections, the
steps sketched above are examined in turn.

\section{Assumptions and preliminaries}\label{assumptions}
\setcounter{equation}{0} 
We assume henceforth that \Ref{finite-jumps} and~\Ref{finite-rates} are satisfied.
Since the index~$j\in\integ_+$ is symbolic in nature, we fix
an $\n\in\rr$, such that~$\n(j)$ reflects in some sense the `size' of~$j$: 
\eq\label{nu-limits}
  \n(j) \ge 1 \ \mbox{for all}\ j\ge0 \quad\mbox{and}\quad
     \lim_{j\to\infty}\n(j) = \infty. 
\en
We then assume that most indices are large and that most transitions involve some
large indices, in the sense that, for $\TT_M := \{j\colon\,\n(j) \le M\}$ and 
$\jj_M := \{J\colon\, J^j=0 \mbox{ for all } j \notin \TT_M \}$,
we have
\eq\label{T+J-assns}
    |\TT_M| \le n_1M^{\b_1};\qquad|\jj_M| \le n_2M^{\b_2},
\en  
for some $n_1,n_2$ and $\b_1 \le \b_2$; note that in fact $\b_2 \le 2\b_1 J_*$ also. 
As a consequence of these assumptions,
for any $s > \b_1$, there exists $K_s < \infty$ such that
\eq\label{nu-sums-1}
  \sjo \n(j)^{-s} \ <\ \infty;\qquad\sum_{j \notin \TT_M} \n(j)^{-s} \ <\ K_s M^{-s+\b_1};
\en
moreover, if $\n(J) := \max_{\{j\colon\,J^j \ne 0\}} \n(j)$, then, for any $s > \b_2$,
\eq\label{nu-sums-2}
  \sjj \n(J)^{-s} \ <\ \infty;\qquad \sum_{J \notin \jj_M} \n(J)^{-s} \ <\ K'_s M^{-s+\b_2},
\en  
for some $K'_s < \infty$. 

\bsk\nin{{\bf Moment assumptions}}

In the proofs that follow, it is important to be able to show that~$x^N$ is largely
concentrated on indices~$j$ with~$\n(j)$ not too large.  This is shown to be the
case in [BL, Section~2], under the following `moment' assumptions.  Defining
$S_r(x) := \sjo x^j \{\n(j)\}^r$, $x\in\rr_0$, and then
\eq\label{UV-defs}
    U_{r}(x) \Def \sjj \a_J(x)\! \Bigl( \sjo J^j \{\n(j)\}^r \!\Bigr); \quad 
    V_{r}(x) \Def \sjj \a_J(x)\! \Bigl( \sjo J^j \{\n(j)\}^r \!\Bigr)^2, \quad  
\en
$x\in\rr$, the assumptions that we need are as follows.  

\begin{assumption}\label{ap-assns}
For~$\nu$ as above, assume that there exist $\rmax\ui, \rmax\ut \ge 1$ 
such that, for all~$X\in\xx_+$, 
\eqa
     &&\sjj \a_J(N^{-1}X) \Blm \sjo J^j \{\n(j)\}^r \Brm \ <\ \infty,\qquad 0\le r\le \rmax\ui,
       \label{3.5a} 
\ena
and also that, for some non-negative constants $k_{rl}$, the inequalities
\eqa   
    U_{0}(x) &\le&  k_{01} S_0(x) + k_{04}, \non \\
    U_{1}(x) &\le&  k_{11} S_1(x) + k_{14},  \qquad\hskip1in 2\le r\le \rmax\ui, \label{U-bnd} \\
    U_{r}(x) &\le&  \{k_{r1} + k_{r2} S_0(x)\} S_r(x) + k_{r4},
    \non
\ena
and
\eq
   \mbox{}\qquad\begin{array}{rl}
    V_{0}(x) &\le\ k_{03}S_1(x) + k_{05}, \\[1ex]
    V_{r}(x) &\le\ k_{r3} S_{p(r)}(x) + k_{r5}, 
   \end{array} \qquad\hskip0.9in 1\le r\le \rmax\ut, \label{V-bnd}
\en
are satisfied,
where $1\le p(r) \le \rmax\ui$ for $1\le r\le \rmax\ut$.
\end{assumption}

As a result of these assumptions, it is shown in [BL, Lemma~2.3 and Theorem~2.4] that, if~$r$
is such that $1\le r\le \rmax\ut$ and if 
$\max\{S_r(x_0^N), S_{p(r)}(x_0^N)\} \le C$ for some~$C$,
then there are constants $C_1$ and~$C_2$, depending on $C,r$ and~$T$, such that
\eq\label{nu-moment-bnd}
   \pr\Bigl[\stT S_r(x_t^N) > C_1 \Bigr] \Le N^{-1}C_2.
\en
 
\msk\nin{{\bf Semigroup assumptions}}

In order to make sense of~\Ref{x-eqn}, we need some assumptions about~$A$.
We assume that
\eq
   A_{ij} \ge 0 \ \mbox{for all}\ i \neq j \ge 0;\qquad
      \sjni A_{ji} < \infty\  \mbox{for all}\ i\ge0, \label{A-assn-1} 
\en
and that, for some $\m \in \R_+^{\Zp}$ such that $\m(m) \ge 1$ for each
$m\ge0$, and for some $w\ge0$,
\eq\label{A-assn-2}
   A^T \m \Le w\m.
\en
We then use~$\m$ to define the $\m$-norm 
\eq\label{Rmu}
    \nm{\x} \Def \smo \m(m)|\x^m| \quad\mbox{on}\quad
    \rr_\m \Def \{\x\in\rr\colon\,\nm{\x} < \infty\},
\en
and, under these assumptions, the transition semigroup~$R$ is well defined
[BL, Section~3], and
\eq\label{R-growth}
     \sio \m(i)R_{ij}(t) \Le \m(j)e^{wt} \qquad \mbox{for all}\ j\ \mbox{and}\ t.
\en
Note that there may be many possible choices for~$\m$, but that we also require 
that~$F$ is locally Lipschitz in the $\m$-norm, in order to ensure that~\Ref{x-eqn}
has a $\m$-continuous solution: we assume that, for any $z>0$,
\eq\label{F-assn-2}
  \sup_{x \neq y\colon\,\nm{x},\nm{y} \le z}\nm{F(x)-F(y)}/\nm{x-y} \Le K(\m,F;z)
     \ <\ \infty,
\en 
and this should be borne in mind when choosing~$\m$. We further assume that, for some 
$\b_3,\b_4$, 
\eq\label{mu-A-growth}
    \m(j) \Le \n(j)^{\b_3} \quad\mbox{and}\quad  |A_{jj}| \Le \n(j)^{\b_4}.  
\en

\bsk\nin{{\bf Transition rate assumptions}}

We need to ensure that the sum of the transition rates,
even when weighted by largish powers of~$\n(j)$, remains bounded.  
To ensure this, we assume that, for some $r_0$ large enough, there exist
$r_1 \le \rmax\ut$, $b\ge1$ and $k_1,k_2>0$ such that
\eq\label{zeta-bnd}
   \sjj \a_J(x)\sjo |J^j|\{\n(j)\}^{r_0} \Le \{k_1 S_{r_1}(x) + k_2\}^b;
\en
this assumption is a specialized version of~[BL, (2.25)].
In view of~\Ref{nu-moment-bnd}, this implies that, if $\max\{S_{r_1}(x_0^N), 
S_{p(r_1)}(x_0^N)\} \le C$, then there are constants $C_1$ and~$C_2$
depending on $C$ and~$T$, such that
\eq\label{moment-rate-bnd}
    \pr\Bigl[\stT \sjj \a_J(x_t^N)\sjo |J^j|\{\n(j)\}^{r_0} > C_1 \Bigr] \Le N^{-1}C_2.
\en
We shall therefore assume that the initial condition needed for~\Ref{moment-rate-bnd}
is indeed satisfied: that, for some~$C<\infty$,
\eq\label{ic}  
   \max\{S_{r_1}(x_0^N), S_{p(r_1)}(x_0^N)\} \le C.
\en
It can be seen from the statement of Theorem~\ref{diffusion-theorem} that the larger we can 
take~$r_0$ in~\Ref{zeta-bnd}, the sharper the approximation bound that we get in~\Ref{CLT-1}, in that
$\z$ can be taken smaller for a given value of the product~$\z r_0$, resulting in larger values
of~$b_1(\z)$.

 Since it is immediate that
\[
   \sjj \sjo  |J^j| \{\n(j)\}^r \a_J(x) \ \ge \sjj \{\n(J)\}^r \a_J(x),
\] 
it follows that, for any $r,s \ge 0$,
\eqs 
    \sum_{J \notin \jj_M} \{\n(J)\}^r \a_J(x) 
      &\le& M^{-s} \sum_{J \notin \jj_M} \{\n(J)\}^{r+s} \a_J(x) \\
      &\le& M^{-s} \sjj \sjo |J^j| \{\n(j)\}^{r+s} \a_J(x),
\ens      
so that, if $r+s \le r_0$,  \Ref{moment-rate-bnd} implies that
\eq\label{out-of-M-bnd-2}
         \stT\sum_{J \notin \jj_M}  \{\n(J)\}^r \a_J(x_t^N) \Le C_1 M^{-s},
\en
except on an event of probability  at most $C_2 N^{-1}$.

\bsk\nin{{\bf Smoothness assumptions}}

We need some smoothness conditions on the rates near the deterministic path~$x$.
First, for some $\d>0$, we assume that~$F$ has second order partial derivatives in the tube
\eq\label{delta-sausage-def}
     B(t,x,\d) \Def \{z \in \xx\colon\, \nm{z-x_s} \le \d \mbox{ for some } 0\le s\le t\},
\en
where~$x$ solves~\Ref{x-eqn}, and that, for any $j,k,l$,
\eq\label{C2-assn}
    \sup_{z \in B(x,t,\d)} \left|{D_{kl} F^j}(z)\right|
     \Le v_{jkl},
\en
where the~$v_{jkl}$ are such that
\eq\label{v-bnd}
     \sjo \m(j) v_{jkl} \Le K_{F2}\m(k)\m(l),
\en
for some $K_{F2} < \infty$. Note that~\Ref{v-bnd} is satisfied if
\eq\label{v-cond1}
  v_{jkl} \Le v^j\m(k)\m(l),\qquad \mbox{where}\quad \nm{v} < \infty.
\en
It is also true under the following condition: that, for each $k$, there exists
$N(k) \subset \Zp$ with $|N(k)| \le n_0$ and $\max_{j\in N(k)}\m(j) \le K_0\m(k)$
such that
\eq\label{v-cond2}
    v_{jkl} \Le K_1\{\m(l)\bone_{\{N(k)\}}(j) + \m(k)\bone_{\{N(l)\}}(j)\};
     \qquad v_{jkl} \Eq 0 \quad\mbox{otherwise},
\en
for suitable $n_0,K_0$ and~$K_1$, all finite.
The first derivative of~$F$ has already been assumed to be $\m$-Lipschitz in~\Ref{F-assn-2};
with the assumption~\Ref{C2-assn}, $F$ becomes continuously $\m$-differentiable in the tube, 
so that, for some constant~$K_{F1}$,
\eq\label{DF-assn}
   \stT \sjo \m(j) |D_k F^j(x_t)| \Le K_{F1}\m(k) \quad\mbox{for all}\quad k.
\en
We also assume that the individual transition rates~$\a_J$ are uniformly $\mu$-Lipschitz 
in~$B(T,x,\d)$, with
\eq\label{alpha-lip}
   \sup_{z_1,z_2\in B(T,x,\d)}|\a_J(z_1)-\a_J(z_2)|/\|z_1-z_2\|_\m \Le K_\a\{\n(J)\}^{\b_5}
\en
for some $K_\a,\b_5 > 0$. This assumption, and those on the second derivatives of~$F$,
go beyond what is required for the law of large numbers in~[BL]; the same is true of
the assumptions \Ref{T+J-assns} and~\Ref{mu-A-growth}.

\bsk\nin{\bf Preliminary conclusions}

We now assume, in addition, that we can take 
\eq\label{r0-assn}
   r_0 \ >\ 2(\b_1+\b_3+\b_4)
\en
in~\Ref{zeta-bnd}.
Then, under the assumptions of this section, it follows from [BL, Theorem~4.7],
with $\z(j) := \{\n(j)\}^{r_0}$, that
the following result holds: for each $T > 0$, 
there exist constants $K_T^{(1)}, K_T^{(2)}$ and $K_T^{(3)}$ such 
that, for all~$N$ large enough, if  
\eq\label{ic-2}
  \|x_0^N - x_0\|_{\mu} \le K_T^{(1)} \sqrt{\frac{\log N}{N}},
\en
then
\begin{equation}
\label{main-result-1}
    \pr \Bigl( \sup_{0 \le t \le T} \|x_t^N - x_t \|_{\mu} 
      > K_T^{(2)} \sqrt{\frac{\log N}{N}}   \Bigr) \le K^{(3)}_T \frac{\log N}{N}.
\end{equation}
We shall from now on also assume that~\Ref{ic} holds with $x$ for~$x^N$. 
Since then~$x$ can be represented as a limit of processes~$x^M$ satisfying~\Ref{ic},
because of~\Ref{main-result-1}, it follows in view of
\Ref{zeta-bnd} and~\Ref{moment-rate-bnd} that we also have
\eq\label{moment-rate-limit-bnd}
    \stT \sjj \a_J(x_t)\sjo |J^j|\{\n(j)\}^{r_0} \Le C_1,
\en 
and therefore, as for~\Ref{out-of-M-bnd-2}, for $r + s \le r_0$,
\eq\label{out-of-M-bnd-1}
       \stT\sum_{J \notin \jj_M}  \{\n(J)\}^r \a_J(x_t) \Le  C_1 M^{-s}.
\en
When approximating~$x^N$ by a deterministic path~$x$, it is natural to choose their
initial values to be close, as in~\Ref{ic-2}.  The impact of also assuming~\Ref{ic} for the 
initial values of both paths is to specify how much closer the components need to be, whose 
indices~$j$ have~$\n(j)$ large.

\bsk\nin{\bf Example}

In the model of Arrigoni~(2003) presented in the introduction, we can take $\n(j) = j+1$,
in which case $\b_1=1$ and~$\b_2 = 2$, the latter because of the migration transition.
Calculation shows that~\Ref{U-bnd} is satisfied for all~$r$, as is~\Ref{V-bnd} also,
with $p(r) = 2r$, so that we can take $\rmax\ui = \rmax\ut = \infty$.  Furthermore,
\Ref{zeta-bnd} is satisfied for any~$r_0$, with~$r_1 = r_0+1$.  The quantities $A$ and~$F$
are given by
\eqs
   A_{ii} &=& -\{\k + i(b_i+d_i+\g)\};\  A^T_{i,i-1} \Eq i(d_i+\g);
     \  A^T_{i,i+1} \Eq i b_i, \quad  i\ge1; \\       
   A_{00} &=& -\k, 
\ens
with all other elements of~$A$ equal to zero, and, writing $s(x) := \sji jx^j$, 
\[
    F^i(x) \Eq \r\g (x^{i-1} - x^i)s(x),\quad i\ge1; \qquad  
    F^0(x) \Eq -\r\g x^0 s(x) + \k,
\]
where we have used the fact that $\sjo x^j = 1$.
Hence Assumption \Ref{A-assn-1} is immediate, and Assumption~\Ref{A-assn-2}
holds for $\m(j) = j+1$ (so that $\b_3=1$), with $w = \max_i(b_i-d_i-\g-\k)_+$ 
(assuming the $b_i$'s and $d_i$'s to be such that this is finite). 
The value of~$\b_4$ depends on the particular choice of the $b_i$ and~$d_i$.
For instance, the stochastic version of Ricker's~(1954) model has both the
$b_i$ and the~$d_i$ uniformly bounded, in which case we can take $\b_4=1$.
However, in the stochastic analogue of Verhulst's~(1838) logistic model, the~$d_i$ 
grow linearly with~$i$, and then one needs $\b_4=2$. 

With the above choice of~$\m$, $F$ can easily be seen to be locally Lipschitz in 
the $\m$-norm, with $K(\m,F;z) \le 4\r\g z$.
The partial derivatives of~$F$ are given by
\eqs
  D_kF^i(x) &=& \r\g k(x^{i-1}-x^i) + \r\g s(x)\{\bone_{\{k\}}(i-1) - \bone_{\{k\}}(i)\};\\
  D_{kl}F^i(x) &=& \r\g\{k[\bone_{\{l\}}(i-1) - \bone_{\{l\}}(i)] 
                + l[\bone_{\{k\}}(i-1) - \bone_{\{k\}}(i)]\},
\ens
for any $i,k,l \ge 0$ (we take $x^{-1}=0$). From this, it follows (using the
elementary bound $j+1 \ge 2j$ in $j\ge0$) that we can take $K_{F1} = 6\r\g
\stT \|x_t\|_\m$ in~\Ref{DF-assn}, and that~\Ref{v-cond2} is satisfied with
$n_0=2$, $K_0=2$ and~$K_1=\r\g$, so that~\Ref{v-bnd} is also satisfied
(one can in fact take $K_{F2} = 4\r\g$).  Finally,
\Ref{alpha-lip} is satisfied, with~$\b_5 = 0$ if the~$d_i$ are uniformly
bounded, and with~$\b_5 = 1$ if they grow linearly, and with~$K_\a$ of the
form $K'(1 + \stT \|x_t\|_\m)$.

\section{Controlling $\h^N$}\label{step-1}
\setcounter{equation}{0}

From now on, we assume that all the assumptions of Section~\ref{assumptions} are
in force.
We first show that the effect of the perturbation~$\h^N$ is negligible.
For~$\h_t^N$, from~\Ref{eta-def}, we need to consider the difference
\[
   \int_0^t R(t-s)\{F(x_s^N) - F(x_s) - DF(x_s)[x_s^N - x_s]\}\,ds.
\]
We note first that, if $\nm{h} \le \d$ for~$\d$ as in Condition~\Ref{delta-sausage-def},
then, from \Ref{v-bnd},
\eqs
    \nm{F(x_s^N) - F(x_s) - DF(x_s)[h]} &\le& \sjo \m(j) \sko\slo |h_kh_l|v_{jkl}
      \Le K_{F2}\nm{h}^2.
\ens
Hence, from~\Ref{main-result-1} and from~\Ref{R-growth}, for all~$N$ 
large enough to ensure that $K_T\ut N^{-1/2}\sqrt{\log N} \le \d$, we have
\eq\label{eta-controlled}
   \sup_{0\le s\le t}\nm{\h_s^N} \Le K_{F2}t\{K_T\ut\}^2e^{wt}N^{-1/2}\log N,
\en
for all $0 < t \le T$, except on a set of probability at most $K_T\uh N^{-1}\log N$.

\section{Discrete to diffusion}\label{step-2}
\setcounter{equation}{0}
We now show that $N^{1/2}m^N$ is close in the $\m$-norm to the diffusion~$W$, given
by
\[
    W_t \Def \sjj JW_J(A_J(t)),
\]
as in~\Ref{W-def}.  We first need to show that this~$W$ indeed has
paths in~$\rr_\m$.  For this, it is enough to show that
\eq\label{W-norm-bnd}
   \sjo \m(j) \sjj |J^j| |W_J(A_J(t))| \ <\ \infty
\en
for all~$t$.

We begin by noting that, using the reflection principle,
if~$B$ is standard Brownian motion, then there exists a
constant~$\g < \infty$ such that, for all $a>0$,
\eq\label{BM-bnd}
   \pr[\sup_{0\le x\le 1} |B(x)| > a\g] \Le e^{1-a^2}.
\en
Thus, from~\Ref{BM-bnd}, for any $C>1$ and $p,T>0$, we have
\[
   |W_J(A_J(t))| \ =_d\  A_J(T)^{1/2} |B(A_J(t)/A_J(T))|
    \Le A_J(T)^{1/2}\g\sqrt{p\log(C\n(J))},
\]
for all $0\le t\le T$, except on a set of probability at most
$e\{C\n(J)\}^{-p}$.  Hence it follows that
\eq\label{W-compt-bnd}
   |W_J(A_J(t))| \Le A_J(T)^{1/2}\g\sqrt{p\log(C\n(J))}
\en
for all $0\le t\le T$ and for all~$J\in\jj$ , except on a set of probability at most
$eC^{-p}\sjj \{\n(J)\}^{-p} < \infty$, by~\Ref{nu-sums-2}, if $p > \b_2$.

For $x,y \ge 0$, one has $\sqrt{x+y} \le (1+\sqrt x)(1+\sqrt y)$.
Substituting from \Ref{W-compt-bnd} into~\Ref{W-norm-bnd} shows
that $\nm{W_t} < \infty$ a.s.\ for all $0\le t\le T$, provided that
\[
 \sjo \m(j) \sjj |J^j| \{1 + \sqrt{\log\n(J)} \} A_J(t)^{1/2}\  <\  \infty,
\]
since~$C$ is arbitrary.
However, by~\Ref{mu-A-growth} and recalling the definition of~$J_*$, we have,
for any $\e>0$,
\eqa
  \lefteqn{ \sjo \m(j)\sjj |J^j| \n(J)^\e A_J(t)^{1/2} }  \non\\
  &\le& \sjj \sjo |J^j| \n(J)^{-\b_{2}/2-\e} \, \n(J)^{\b_3 + 2\e + \b_2/2}  A_J(t)^{1/2} 
       \phantom{HHHHHHHHHHHH}\non\\
  &\le& J_*\Blb\! \Bl\sjj \n(J)^{-\b_{2}-2\e}\Br^{1/2}
     \Bl\sjj\int_0^T  \n(J)^{\b_2+2\b_3+4\e} \a_J(x_s)\,ds \Br^{1/2}\! \Brb,\label{sqrt-calc}    
\ena
and both sums in the final expression are finite, by~\Ref{nu-sums-2} and~\Ref{moment-rate-limit-bnd},
provided that $\b_2+2\b_3 < r_0$ and that~$\e$ is small enough.

Having established that~$W$ indeed has paths in~$\rr_\m$, we now need to show that it is
close to~$N^{1/2}m^N$ in the $\m$-norm, if the Brownian motions~$W_J$ are suitably
chosen.   The relationship between $N^{1/2}m^N$ and~$W$ arises because $N^{1/2}m^N$
can be represented in the form
\eqa\label{m^N-Poisson}
   Nm_t^N &:=& N\Blb x_t^N - x_0^N - \int_0^t \sjj \a_J(x_s^N)\,ds\Brb \non\\
    &=& \sjj J \{P_J(NA_J^N(t)) - NA_J^N(t)\},
\ena
where $A_J^N(t) := \int_0^t \a_J(x_s^N)\,ds$, and the~$P_J$'s are independent Poisson
processes.  Now $\{N^{-1/2}(P_J(Nt) - Nt),\,t\ge0\}$ can be well approximated by
a Brownian motion, and $A_J^N(t)$ is close to $A_J(t) := \int_0^t  \a_J(x_s)\,ds$,
by~\Ref{alpha-lip} and~\Ref{main-result-1}.

We thus wish to show that the~$W_J$ can be chosen in such a way that
\eq\label{mu-norm-diff}
   \stT \sjo \m(j) \Blm \sjj J^j Z_J^N(A_J^N(t)) - \sjj J^j W_J(A_J(t)) \Brm
\en
is small, where we define
\eq\label{Z-def}
    Z_J^N(s) \Def N^{-1/2}\{P_J(Ns) - Ns\}.
\en 
 For use in the next section, we prove somewhat more: that, under
appropriate conditions, we can replace $\m(j)$ in~\Ref{mu-norm-diff} by the
larger quantity $\n_*(j) := \{\n(j)\}^{\b_3+\b_4}$, and still obtain something that is small.
To do so, we begin by bounding the sum by $T_1(t) + T_2(t) + T_3(t)$, where
\eqa
   T_1(t)&:=&  \sum_{J \in \jj_M} \sjo |J^j| \,\n_*(j)\, |Z_J^N(A_J^N(t)) - W_J(A_J(t))|; \non\\
   T_2(t)&:=&  \sum_{J \notin \jj_M} \sjo |J^j| \,\n_*(j)\, 
                                  |Z_J^N(A_J^N(t))|;       \label{main-estimate} \\
   T_3(t)&:=&  \sum_{J \notin \jj_M} \sjo |J^j| \,\n_*(j)\,|W_J(A_J(t))|. \non 
\ena      
Here, $M$ is to be chosen later as~$N^\z$, for some suitable small $\z > 0$.

We begin with~$T_2(t)$, which we deal with by showing that, for suitable choice of~$M$,
$N \sum_{J \notin \jj_M} A_J^N(T)$ is small.  Indeed,
\[
   N \sum_{J \notin \jj_M} A_J^N(T) \Eq N\int_0^T  \sum_{J \notin \jj_M} \a_J(x_u^N)\,du,
\]
which is bounded by using~\Ref{out-of-M-bnd-2} with $r=0$ and $s \le r_0$, together 
with~\Ref{moment-rate-bnd}; the quantity is of order $NM^{-s}$ for any $s \le r_0$,
except on an event of probability of order~$O(N^{-1})$.  Thus, except on an event
with probability of order $O(N^{-1} + NM^{-r_0})$, $N^{-1/2}P_J(NA_J^N(T)) = 0$ for all
$J \notin \jj_M$.  Furthermore, the contribution from the compensators is bounded by
\eqa
    N^{1/2} \sum_{J \notin \jj_M} \n_*(j) |J^j| A_J^N(T) &\le&
        N^{1/2}T\stT \sum_{J \notin \jj_M} \n(j)^{\b_3+\b_4} |J^j| \a_J(x_t^N) \non\\
          &=& O(N^{1/2}M^{-s'}),  \label{compensator-bnd}
\ena
by~\Ref{out-of-M-bnd-2}, if $\b_3+\b_4 + s' \le r_0$, except on an event with probability
of order~$O(N^{-1})$. Recalling~\Ref{r0-assn}, this proves that, for 
any~$s' \le r_0 - \b_3 - \b_4$,
\eq\label{T_2}
    \stT T_2(t) = O(N^{1/2}M^{-s'}),
\en 
except on an event of probability of order $O(N^{-1} + NM^{-r_0})$.
  
For~$T_3(t)$, we use~\Ref{BM-bnd} to give
\eq\label{W_J-bnd}
  \pr\Bigl[\stT |W_J(A_J(t))| > \{A_J(T)\}^{1/2}\g\sqrt{p\log\n(J)}\Bigr] \Le e\{\n(J)\}^{-p},
\en  
for any $p>0$.  Hence, for any $p > \b_2$, it follows that
$|W_J(A_J(t))| \le \{A_J(T)\}^{1/2}\g\sqrt{p\log\n(J)}$ for all $J\notin \jj_M$ and for all 
$0\le t\le T$,
except on an event~$E_3$ of probability of order $O(M^{-p+\b_2})$, from~\Ref{nu-sums-2}.
But then, except on~$E_3$, for all $0\le t \le T$,
\eqs
   \sum_{J \notin \jj_M} \sjo |J^j| \n_*(j) |W_J(A_J(t))| 
       &\le& \sum_{J \notin \jj_M} \sjo |J^j| \n_*(j)  \{A_J(T)\}^{1/2}\g\sqrt{p\log\n(J)} \\     
       &\le& K_\e \sqrt p \sum_{J \notin \jj_M} J_* \n(J)^{\b_3+\b_4+\e}  \{A_J(T)\}^{1/2},
\ens
for any $\e > 0$, with suitable choice of~$K_\e$.  But now, for any $r > 0$,
\eqa
   \lefteqn{\sum_{J \notin \jj_M}  \n(J)^{\b_3+\b_4+\e}  \{A_J(T)\}^{1/2}} \non\\ 
       &&\Le \Blb  \sum_{J \notin \jj_M}  \n(J)^{-2(r-\b_3-\b_4-\e)} \Brb^{1/2}
              \Blb  \sum_{J \notin \jj_M}  \n(J)^{2r} A_J(T) \Brb^{1/2} \phantom{XXX}\non \\[2ex]
      &&\Le \{K'_{2(r-\b_3-\b_4-\e)}\}^{1/2}
          \Blb   M^{-2(r-\b_3-\b_4-\e)+\b_2} \Brb^{1/2} \Blb T  M^{-2r'} \Brb^{1/2},
       \label{sqrt-sum}  
\ena
by \Ref{nu-sums-2} and~\Ref{out-of-M-bnd-1}, so long as 
$r > \b_3+\b_4+\e + \b_2/2$ and $2(r+r') \le r_0$.
Hence, if $r_0 > \b_2 + 2(\b_3+\b_4)$, then for any $r' <  (r_0 - \b_2)/2 - (\b_3+\b_4)$ we have 
\eq\label{T_3}
   p^{-1/2}\stT T_3(t) = O(M^{-r'}),
\en  
with an implied constant uniform for all $p > \b_2$,
except on an event with probability of order $O(M^{-p+\b_2})$.

So far, the bounds have been achieved without any specific choice of the Brownian 
motions~$W_J$, but, for~$T_1(t)$, we need to be more precise.  We treat each~$J$ 
separately, since the underlying
Poisson processes~$P_J$ are independent, and match the centred and normalized
Poisson process $Z_J^N$ to a Brownian motion~$W_J$ using the KMT construction.
We need only to do this over a limited time interval, since, from~\Ref{moment-rate-bnd},
\[
    \stT  \sjj   \a_J(x_t^N) \Le C_{1},
\]
except on an event~$E_0$ of probability of order~$O(N^{-1})$, so that, off~$E_0$,
\eq\label{A_J^N-bnd}
   A_J^N(T) \Le TC_{1} \quad \mbox{for all} \ J.
\en
We use Koml\'os, Major \& Tusn\'ady~(1975, Theorem~1 (ii)), together with~\Ref{BM-bnd}
to interpolate between integer time points, applied to the centred unit rate Poisson
process. This implies that, for any $p>0$, we can choose~$W_J$ in such a way that
\[
   \sup_{0 \le t \le TC_{1}} \{N^{-1/2}(P_J(Nt) - Nt) - W_J(t)\} \Le k_p N^{-1/2}\log N,
\]
for a constant~$k_p$,
except on an event~$\tE_{Jp}$ of probability of order~$O(N^{-p})$.  Thus the same bound
holds for all $J\in\jj_M$ except on an event~$\tE_p$ of probability of order
$O(M^{\b_2}N^{-p})$.  Hence, except on $E_0 \cup \tE_p$, an event of
probability of order $O(N^{-1} + M^{\b_2}N^{-p})$, we have 
\eqa
    T_{11}(t) &:=& \sum_{J \in \jj_M} \sjo |J^j| \,\n_*(j)\, |Z_J^N(A_J^N(t)) - W_J(A_J^N(t))| \non \\
    &\le& J_*\sum_{J \in \jj_M} \{\n(J)\}^{\b_3+\b_4} \,k_p N^{-1/2}\log N \non\\
    &=& O(M^{\b_2+\b_3+\b_4}N^{-1/2}\log N), \label{T_{11}}
\ena
for all $0\le t\le T$.  It thus remains to bound
\eq\label{T_{12}-def}
   T_{12}(t) \Def \sum_{J \in \jj_M} \sjo |J^j| \,\n_*(j)\, |W_J(A_J^N(t)) - W_J(A_J(t))|.
\en

First, note that, by~\Ref{alpha-lip},
\eq\label{AJ-diff}
   \stT |A_J^N(t) - A_J(t)| \Le TK_\a\{\n(J)\}^{\b_5}\stT\nm{x_t^N - x_t},
\en
and that 
\eq\label{x-norm-bnd}
    \stT\nm{x_t^N - x_t} \Le K_T^{(2)}N^{-1/2}\sqrt{\log N},
\en
by~\Ref{main-result-1},
except on an event with probability of order~$O(N^{-1}\log N)$.  Furthermore, 
by~\Ref{moment-rate-limit-bnd}, $A_J(T) \le TC_{1}$ for all~$J$.  Then, for a Brownian
motion~$W$ and for $0 < \d < 1$, by a standard argument based on~\Ref{BM-bnd},
\eq\label{W_J-path-bnd}
   \pr\left[ \sup_{0\le u \le A, |s|\le\d} |W(u+s) - W(u)| > 3\g\d^{1/2}\sqrt{r\log(1/\d)} \right]
    \Le (A\d^{-1}+2)e\d^r,
\en
for any~$r>0$, with $\g$ chosen as for~\Ref{BM-bnd}.  Hence, taking $W = W_J$ and, in
view of \Ref{AJ-diff} and~\Ref{x-norm-bnd}, taking
$$
   \d \Eq \d_J \Eq TK_\a K_T\ut \{\n(J)\}^{\b_5} N^{-1/2}\sqrt{\log N}
$$ 
and $A=TC_{1}$
in~\Ref{W_J-path-bnd}, it follows that, for any $\e,r'>0$, there is a $K_{r'\e} < \infty$ such that
\[
    \stT |W_J(A_J^N(t)) - W_J(t)| \Le K_{r'\e}M^{\b_5/2} N^{-1/4+\e} \quad\mbox{for all}\quad J\in\jj_M,
\]
except on an event of probability of order~$O(N^{-1}\log N + M^{\b_2+2\b_5}\{M^{2\b_5}/N\}^{r'})$.  
Off the exceptional event, we have
\eqa
   \stT T_{12}(t) &\le&  J_* \sum_{J \in \jj_M}  \n_*(J) K_{r'\e}M^{\b_5/2} N^{-1/4+\e}   \non\\
    &=&  O(M^{\b_2+\b_3+\b_4+\b_5/2}N^{-1/4+\e}), \label{T_{12}}
\ena
for any $\e > 0$, and the exceptional event can be made to have probability of order~$O(N^{-1}\log N)$
by choosing~$r'$ large enough, provided that~$M$ is bounded by a small enough power of~$N$.

Combining \Ref{T_2}, \Ref{T_3}, \Ref{T_{11}} and~\Ref{T_{12}},  and choosing $M = N^\z$,
we see that we have no useful bound unless $\z r_0 > 1$ (because of the exceptional
event in~\Ref{T_2})
and $\z < 1/\{4(\b_2+\b_3+\b_4)+2\b_5\}$ (in view of~\Ref{T_{12}}), so
that $r_0 > 4(\b_2+\b_3+\b_4)+2\b_5$ is a minimal requirement.  Note
that this assumption on~$r_0$ is more restrictive than that assumed in 
Section~\ref{assumptions}. The error bound in~\Ref{T_{11}} is always smaller than 
that in~\Ref{T_{12}}, and with $\z r_0 > 1$, the error bound in~\Ref{T_2} is smaller 
than that in~\Ref{T_3}.   This translates into the following conclusion:
if $r_0 > 4(\b_2+\b_3+\b_4)+2\b_5$, then for any $1/r_0 < \z < 1/\{4(\b_2+\b_3+\b_4)+2\b_5\}$  
we have
\eqa
  \lefteqn{\stT\nm{N^{1/2}m_t^N - W_t}} \non\\ &&\Le
   \stT \sjo \{\n(j)\}^{\b_3+\b_4} \Blm \sjj J^j Z_J^N(A_J^N(t)) - \sjj J^j W_J(A_J(t)) \Brm \non\\
       &&\Eq O(N^{-b_1})  \label{mu-norm-diff-2}
\ena
except on an event of probability of order~$O(N^{-b_2})$, for any
\eqa
    b_1 &<& b_1(\z) \phantom{HHHHHHHHHHHHHHHHHHHHHHHHHHHHH}\non\\
   &:=& 
      \min\{\halfh - \z(\b_2+\b_3+\b_4+\b_5/2), \half\z(r_0 - \b_2 - 2(\b_3+\b_4))\}; 
         \label{exponents} \\  
    b_2 &<& b_2(\z) \Def \min\{\z r_0-1,1\}.\non
\ena

\section{The existence of~$\tW$, and its approximation}\label{W-tilde}
\setcounter{equation}{0}
The next step in the argument is to show that the process~$\tW$ in~\Ref{tW-def}, 
related to~$W$ exactly as~$\tm^N$ is related to~$m^N$ through~\Ref{mN-tilde-def}, 
is well defined, and that it is indeed the limiting analogue of the process 
$N^{1/2}\tm_t^N$.  For its existence, recalling~\Ref{tW-def}, it is enough to 
show that $R(t-s)A W_s$ 
exists for each $s,t$, and belongs to~$\rr_\m$.  For this, it is enough to show
that
\[
   \sio \m(i) \sjo R_{ij}(t-s) \sko |A_{jk}| \sjj J_k|W_J(A_J(s))|
\]
is a.s.\ bounded.  Now, in view of \Ref{A-assn-2}, \Ref{R-growth} and~\Ref{mu-A-growth} and recalling 
that the off-diagonal entries of~$A$ are non-negative, this will be the case if we 
can bound
\eq\label{bnd-1}
   \sko \m(k)(1+|A_{kk}|) \sjj |J_k| |W_J(A_J(s))| 
       \Le J_* \sjj \{\n(J)\}^{\b_3+\b_4}|W_J(A_J(s))|
\en
uniformly in~$s$.  Now, once again from~\Ref{BM-bnd}, for any $C,r>0$,
\[
   \pr\left[\stT |W_J(A_J(t))| > \{A_J(T)\}^{1/2}\g \sqrt{r\log(C\n(J))} \right]
     \Le e \{C\n(J)\}^{-r}, 
\]
so that, in view of~\Ref{nu-sums-2}, taking any $r > \b_2$, there is a (random)~$C$ such that 
\[
    \stT |W_J(A_J(t))| \Le \g \sqrt r \{A_J(T)\}^{1/2}\sqrt{\log(C\n(J))}
\]
a.s.\ for all~$J$.  But now, returning to~\Ref{bnd-1}, we just need to show that
the quantity
$$
   \sjj \{\n(J)\}^{\b_3+\b_4+\e} \{A_J(T)\}^{1/2}
$$ 
is finite
for some $\e>0$, and this is achieved as in~\Ref{sqrt-sum}, if $r_0 > \b_2 + 2(\b_3+\b_4)$.

To show that $\tW$ is a good approximation to $N^{1/2}\tm^N$, we begin with
the result proved in~\Ref{mu-norm-diff-2} above, that, except on an event of 
probability of order~$O(N^{-b_2})$, $\stT \nm{N^{1/2}m_t^N - W_t} = O(N^{-b_1})$
for any $b_1 < b_1(\z), b_2 < b_2(\z)$.  This quantity is one element 
of~$\nm{N^{1/2}\tm_t^N - \tW_t}$; the other is
\[
    \int_0^t R(t-s) A(N^{1/2}m_s^N - W_s)\,ds.
\]
Arguing much as in the previous paragraph, we need to bound
$$
    \stT\sjj \{\n(J)\}^{\b_3+\b_4}|W_J(A_J(t)) - Z_J^N(A_J^N(t))|.
$$
But this is exactly what we achieved in~\Ref{mu-norm-diff-2}. Hence, for~$\z$ such that
$1/r_0 < \z < 1/\{4(\b_2+\b_3+\b_4)+2\b_5\}$ and for any $b_1 < b_1(\z)$, $b_2 < b_2(\z)$,
\eq\label{W-tilde-appx}
   \stT\nm{\tW_t - N^{1/2}\tm_t^N}  \Eq  O(N^{-b_1}),
\en
except on an event of probability of order $O(N^{-b_2})$.

\section{The final approximation}\label{final}
\setcounter{equation}{0}
The final step in the argument is to compare the solution~$U^N$ to~\Ref{main-eqn}
with the solution~$Y$ to~\Ref{Y-eqn}.  Both satisfy the general equation
\eq\label{Z-eqn}
   Z_t \Eq R(t)Z_0 + \int_0^t R(t-s)DF(x_s)[Z_s] + z_t,
\en
but with different initial conditions~$Z_0$ and forcing functions~$z$; and their
difference~$Y-U^N$ also satisfies~\Ref{Z-eqn}, with initial conditions and forcing
functions subtracted. Now,
for~$U^N$, the forcing function $\h^N + N^{1/2}\tm^N$ is close to the forcing
function~$\tW$ for~$Y$, because of \Ref{eta-controlled} and~\Ref{W-tilde-appx}, 
and we shall assume that $U^N_0$ and~$Y_0$ are also close to one another,
so that both differences are small.  We now show that this implies that the
difference between $Y$ and~$U^N$ is also small.
 
First, the assumption~\Ref{DF-assn} implies that, for $w\in\rr_\m$, 
$\nm{DF(x_s)[w]} \le K_{F1}\nm{w}$. It is then immediate from~\Ref{R-growth} that
\eqs
   \left\|\int_0^t R(t-s) DF(x_s)[Z_s]\,ds \right\|_\m &\le&
     K_{F1}\int_0^t e^{w(t-s)}\nm{Z_s}\,ds \\
    &\le& K_{F1}e^{wt}\int_0^t \nm{Z_s}\,ds,
\ens
and that $\nm{R(t)Z_0} \le \nm{Z_0}e^{wt}$.  Hence, for $0\le t\le T$, it follows that
\eq\label{Gronwall-bnd}
   \nm{Z_t} \Le \Bigl\{\nm{Z_0}e^{wT} + \sup_{0\le s\le T} \nm{z_s}\Bigr\}e^{Ct},
\en
with $C = K_{F1}e^{wT}$.  Applying~\Ref{Gronwall-bnd} to $Z=Y-U^N$, 
and using the bounds in~\Ref{eta-controlled} and
\Ref{W-tilde-appx}, it follows that, except on an event of probability of order $O(N^{-b_2})$,
\[
    \stT\nm{Y_t-U_t^N} \Eq O(N^{-b_1}),
\]
where $U^N := N^{1/2}(x^N-x)$ and~$Y$ is the solution to~\Ref{Y-eqn}, provided that
$\nm{Y_0-U^N_0} = O(N^{-b_1})$ also.  This proves
the main theorem of the paper:

\ignore{
show that the solution~$Y$ to~\Ref{Y-eqn},
with $\tW$ replaced by an arbitrary path~$y$ in $\rr_\m$, is $\mu$-Lipschitz
in~$y$.  However, this follows directly from the assumption~\Ref{DF-assn},
which implies that, for $z\in\rr_\m$, $\nm{DF(x_s)[z]} \le K_{F1}\nm{z}$.
It is then immediate from~\Ref{R-growth} that
\eqs
   \left\|\int_0^t R(t-s) DF(x_s)[Y_s]\,ds \right\|_\m &\le&
     K_{F1}\int_0^t e^{w(t-s)}\nm{Y_s}\,ds \\
    &\le& K_{F1}e^{wt}\int_0^t \nm{Y_s}\,ds,
\ens
and that $\nm{R(t)Y_0} \le \nm{Y_0}e^{wt}$, so that, for $0\le t\le T$,
\[
   \nm{Y_t} \Le \Bigl\{\nm{Y_0}e^{wT} + \sup_{0\le s\le T} \nm{y_s}\Bigr\}e^{Ct},
\]
with $C = K_{F1}e^{wT}$.  In view of the bounds in~\Ref{eta-controlled} and
\Ref{W-tilde-appx}, it follows that, except on an event of probability of order $O(N^{-b_2})$,
\[
    \stT\nm{Y_t-U_t^N} \Eq O(N^{-b_1}),
\]
where $U^N := N^{1/2}(x^N-x)$ and~$Y$ is the solution to~\Ref{Y-eqn}.  This proves
the main theorem of the paper:
}

\begin{theorem}\label{diffusion-theorem}
Suppose that the assumptions of Section~\ref{assumptions} are satisfied, and that we
can take $r_0 > 4(\b_2+\b_3+\b_4)+2\b_5$ in~\Ref{zeta-bnd}.  
For any~$\z$ such that $1/r_0 < \z < 1/\{4(\b_2+\b_3+\b_4)+2\b_5\}$,
define 
\[
   b_1(\z) \Def \min\{\halfh - \z(\b_2+\b_3+\b_4+\b_5/2), \half\z(r_0 - \b_2 - 2(\b_3+\b_4))\},
\]   
and $b_2(\z) \Def \min\{\z r_0-1,1\}$.
Suppose that $\nm{Y_0 - U^N_0} \le KN^{-b_1(\z)}$.
Then, for any $b_1 < b_1(\z)$, $b_2 < b_2(\z)$,
we can construct copies of $Y$ and~$U^N$ on the same probability space, in such a way
that
\[
    \stT\nm{Y_t-U_t^N} \Eq O(N^{-b_1}),
\]
except on an event whose probability is of order~$O(N^{-b_2})$. 
\end{theorem}

So, for example, in the model of Arrigoni~(2003), we can take~$r_0$ as big as we wish,
and then~$\z r_0 = 2$, allowing~$b_1 = 1/4 - \e$  and $b_2 = 1 - \e$ for any $\e>0$.
However, these rates can only be attained for correspondingly well controlled
initial conditions: in addition to~\Ref{ic-2}, it is necessary to ensure that
\Ref{ic} is satisfied, so that $S_{2(r_0+1)}(x_0^N) \le C$ for some~$C<\infty$ 
and for all~$N$,  and that $S_{2(r_0+1)}(x_0) \le C$ also.  For stochastic logistic
dynamics within the patches, with $b_i = b$ and $d_i = d + ci$, we need to take
$r_0$ to exceed $4(\b_2+\b_3+\b_4)+2\b_5 = 22$ to yield an error bound that
converges to zero with~$N$, and thus require the initial conditions to have 
uniformly bounded $(46+\d)$-th moments for some $\d>0$.

\section*{Acknowledgement}
The authors wish to thank the Institute for Mathematical Sciences of the
National University of Singapore and the University of Melbourne
for providing welcoming environments while part of this work was
accomplished.  MJL also thanks the University of Z\"urich and ADB Monash University
for their hospitality on a number of visits.

\end{document}